\documentclass[11pt]{amsart}
\usepackage{geometry}                
\usepackage[parfill]{parskip}    
\usepackage{graphicx}
\usepackage{amssymb}
\usepackage{epstopdf}
\usepackage[all, 2cell]{xy}
\usepackage{enumitem}
\usepackage[colorlinks]{hyperref}
\DeclareRobustCommand{\subtitle}[1]{\\#1}

\title{One Mathematic(s) or Many? \subtitle{Foundations of Mathematics in Today's Mathematical Practice}}

\author{Andrei Rodin}

\date{\today}

\begin{document}
\maketitle

\renewcommand{\subtitle}[1]{}

\begin{abstract}
The received Hilbert-style axiomatic foundations of mathematics has been designed by Hilbert and his followers as a  tool for meta-theoretical research. Foundations of mathematics of this type fail to satisfactory perform more basic and more practically-oriented functions of theoretical foundations such as verification of mathematical constructions and proofs. Using alternative foundations of mathematics such as the Univalent Foundations is compatible with using the received set-theoretic foundations for meta-mathematical purposes provided the two foundations are mutually interpretable. Changes in foundations of mathematics do not, generally, disqualify mathematical theories based on older foundations but allow for reconstruction of these theories on new foundations.  Mathematics is one but its foundations are many.     
\end{abstract}

\section{``Practical'' Foundations of Mathematics and Foundations of Mathematics \emph{simpliciter}}
Foundations of a theoretical discipline are traditionally conceived as a join between this discipline and its philosophy; accordingly, the 
Philosophy of Mathematics typically focuses on the foundations of mathematics (FOM) rather than on other areas of mathematical research. The Philosophy of Mathematical Practice attempts, among other things, to enlarge the scope of philosophical analysis of mathematics by paying more attention to past and present mathematical developments beyond FOM and, more specifically, beyond the \emph{Grundlagenkrise} that occurred in mathematics at the turn of the 19th and 20th centuries and the related subsequent developments \cite{Mancosu:2008}. Without trying to diminish the importance and philosophical relevance of mathematical practices and contents beyond FOM, I would like to stress that the mismatch between theory and practice in mathematics, which provides a rationale for distinguishing the Philosophy of Mathematical \emph{Practice} into a separate area of research, concerns FOM in the first place.   

Described in few words, the current situation in FOM is as follows. The mainstream FOM is a relatively isolated and very specialised area of today's mathematical research. The majority of working mathematicians and mathematical educators know very little about FOM and cannot see its relevance to their works. This is certainly not how FOM, traditionally conceived, was expected to function. This unusual situation  produces a range of different reactions in the mathematical community. 

Some mathematicians including many  FOM researchers explain and justify the gap between FOM and the rest of mathematics by a natural division of labor. They assume  (i) that an ordinary  working mathematician does not need to master details of the foundations of her discipline and (ii) that the research in FOMü just like research in any other special area of mathematics, requires a long professional training that cannot  possibly be provided for all working mathematicians. A version of this view may go along with the conviction that FOM is of little or no significance for other mathematical disciplines. This view, however, is also compatible with the belief that all mathematical disciplines crucially depend on FOM whether or not working mathematicians are aware of it. 

A number of  mathematicians working outside the mainstream FOM research believe that mathematics needs foundations in a more traditional sense of the term, according to which FOM should represent some core mathematical contents and be known to and used by all mathematicians. On this basis Lawvere and Rosebrugh boldly disqualify the received FOM as a valid foundation:    

\begin{quote}
A foundation makes explicit the essential general features, ingredients, and operations of a science, as well as its origins and general laws of development. The purpose of making these explicit is to provide a guide to the learning, use, and further development of the science. A ``pure" foundation that forgets this purpose and pursues a speculative ``foundations" for its own sake is clearly a nonfoundation.
(\cite{Lawvere&Rosebrugh:2003}, p.235) 
\end{quote}

Some other mathematicians try to make a peaceful agreement between the two concepts of FOM. They reserve the name of \emph{practical foundations}  for the traditional notion without trying to challenge the received FOM explicitly  \cite{Taylor:1999}. On a parallel \textemdash\ albeit chronologically more recent \textemdash\  development in the philosophical camp a group of philosophers  constituted the \emph{Philosophy of Mathematical Practice} as a subfield of philosophy of mathematics outside the mainstream FOM-oriented philosophy of mathematics \cite{Mancosu:2008}. Since pure mathematics is a theoretical discipline \emph{par excellence}, the grounds of the relevant distinction between mathematical theory and mathematical practice are far from clear, moreover so that philosophers of mathematical practice in their work systematically refer to mathematical theories and other mathematical content, but not only to practical issues related to mathematics research and education. The popular wisdom according to which ``in theory there is no difference between theory and practice, while in practice, there is'' perfectly applies in this case. Practically speaking, it is indeed the case that philosophers who associate themselves with the mainstream FOM research and philosophers of mathematical practice are doing different research and form different communities. But from a theoretical point of view it is not clear to me why these lines of philosophical inquiry into mathematics should be so different.

In order to clarify the situation I will not go deeply into philosophical generalities about theory and practice. Instead, I will analyse the nature and the genesis of the aforementioned gap between the received FOM and the rest of mathematics. I will use the above concise non-orthodox description of FOM coined by Lawvere and Rosenburgh as my working definition of this concept.

\section{The Received FOM and Meta-Mathematics}
Today's received concept of FOM stems from Hilbert's seminal works in axiomatic geometry \cite{Hilbert:1899} and his later works on what is known today under the name of \emph{Hilbert's Program}, which includes a modern form of Proof Theory \cite[vol. 2]{Hilbert&Bernays:1934-1939}. A special feature of this conception of FOM, which distinguishes it from more traditional conceptions, is its \emph{meta-mathematical} character. The idea of meta-mathematics, as Hilbert and Bernays explain it in \cite[ch.2, sect. 5]{Hilbert&Bernays:1934-1939}, is to provide \emph{mathematical} answers to certain general mathematically-related questions, which have been earlier seen as purely logical and philosophical \cite{Ferreiros:2005} 
\footnote{Such as Hilbert's \emph{Entscheidungsproblem}.}. 

The idea of meta-mathematics can be traced back to Descartes and Leibniz but it is not constitutive in the traditional FOMs. In particular, Euclid's \emph{Elements}  (which qualify as FOM in our sense) has no meta-mathematical capacity. Let me use this classical example for explaining were the idea of empowering FOM with meta-mathematics comes from. 

Euclid's \emph{Elements} comprise a number of impossibility results such as a proof of impossibility to find a common measure for a pair of arbitrary straight segments\footnote{See the \emph{Elements}, Book 10,  Problem 10.}.  Such proofs show that certain geometrical constructions can be specified in Euclid's theoretical language but, on pain of contradiction, cannot be performed within Euclid's theory. Such impossibility proofs solve certain mathematical problems in the \emph{negative} way by showing that certain desired constructions are impossible. 

However there are also well-known mathematical problems that are are well-posed in Euclid's mathematics but cannot be solved using only Euclid's principles 
either positively (by performing the required construction), or negatively (by showing that the required construction is impossible). The Angle Trisection problem is an example. A trisected angle is not an impossible geometrical construction. One can easily perform an instance of such construction by tripling  a given angle (as usual by ruler and compass). Yet, as we know it today thanks to Galois theory, the trisecting of arbitrary given angle by ruler and compass is not possible. This latter result is meta-mathematical (relatively to Euclid's mathematics) because it belongs to a different theory, to wit Galois theory, built on different principles. In what follows I call a mathematical theory $M_{E}$  \emph{meta}-mathematical with respect to mathematical theory $E$ just in case theorems of $M_{E}$ describe relevant general features of $E$ such as its consistency, limits of its methods, etc. \footnote{This informal description of the concept of meta-mathematics is broader and less exact than formal conceptions found in standard textbooks \cite{Kleene:2009} but it better elucidates the epistemic role of meta-mathematics than more technical definitions.}.

Let us call a theory \emph{epistemically complete} when it allows one to solve \textemdash\ positively or negatively \textemdash\ any well-posed problem within this very theory (see \cite{Pastin:1979} for a formal treatment). The concept of epistemic completeness raises a bunch of questions such as:  Are epistemically \emph{in}complete theories (epistemically) acceptable?; Is the epistemic completeness epistemically valuable?; How does epistemic completeness relate to the semantic completeness and other forms of theoretical completeness? and, last but not least, Are epistemically complete theories possible? David Hilbert was occupied with issues of epistemic completeness during all of his long-term work in FOM. His view on this matter was formed in the context of the so-called  \emph{Ignoramibusstreit} (\emph{ignoramibus} debate; the debate about unknowable) that was started  back in 1872  by Emil du Bois-Reymond who claimed that certain scientific questions and problems are unsolvable in principle. Responding to du Bois-Reymond Hilbert developed an opposite view according to which every well-established theory (and every mathematical theory at the first place) must be epistemically complete. Hilbert's belief remained remarkably stable during his long academic career. Hilbert expressed it publicly on many occasions, most famously in his 1900 address delivered in the Sorbonne at the International Congress of Mathematicians where he  pronounced his famous ``no ignoramibus'' (in mathematics), and also in his last major public lecture delivered in 1930 in K\"onigsberg, which Hilbert concluded with the words ``We must know. We will know'' \cite{McCarty:2005}.

Hilbert believed that the epistemic completeness of mathematics (appropriately reformed) can, and should, be irrefutably proved by rigorous mathematical means, and not only be supported by general philosophical arguments. Whence the idea to identify FOM with a meta-theory of all mathematics, which could rigorously proof the epistemic validity of every mathematical theory deserving the name. Hilbert's meta-mathematical conception of  FOM was based on the assumption that the wanted foundational theory could be kept minimal and form a reliable, trusted core of contemporary mathematics without imposing restrictions onto its conceptual contents. In 1927 Hilbert  speculated that since ``a formalised proof, like a numeral, is a concrete and surveyable object''  proving the formal consistency of a given formalised theory ``is a task that fundamentally lies within the province of intuition just as much as does in contentual number theory the task, say, of proving the irrationality of $\sqrt{2}$'' \cite[p. 471]{Hilbert:1927}. Retrospectively, it is clear that at this point Hilbert seriously underestimated the complexity of metamathematical tasks related to his foundational goals. 

The question of how much (if any) of Hilbert's aspiration remains today sound in view of later FOM-related meta-mathematical and logical research, has been widely discussed since 1930s  \cite{Zach:2003},\cite{Detlefsen:1990},\cite{Artemov:2020}. Leaving aside this and other special questions about the Hilbert-style FOM I would like to stress once again that the very idea Hilbert had of providing FOM with a metamathematical capacity is highly original, and should not be simply taken for granted. From the traditional perspective on FOM, which has been outlined above, a meta-mathematical capacity of FOM appears as useful additional feature. It goes without saying that this new feature of FOM cannot serve its intended epistemic purpose unless the given FOM also performs its more basic functions,  which Lawvere and Rosenburgh describe as ``making explicit the essential general features, ingredients, and operations'' of mathematics. Unless a given FOM performs these basic functions satisfactorily, its meta-mathematical power remains useless and irrelevant.

\section{Is the received FOM adequate?} 
   
Any FOM has two sides. One is normative, logical and epistemological. It deals with questions like What is a correct mathematical proof? and provides tools for building sound mathematical theories and means for ruling out erroneous theories. The other side is ``practical'' and descriptive. It has to do with applications of the given FOM in mathematical practice including mathematics education and research, and, in particular, with the practical enforcement of epistemic norms associated with this given FOM. 

An epistemological analysis of FOM needs to take  both these sides into account. On the one hand, stripping away the normative dimension would reduce the epistemological  analysis to historical, sociological and psychological accounts and thus destroy its specific character and purpose. On the other hand, ignoring the past and present mathematical practice deprives the epistemological analysis of mathematics of its subject-matter. This concerns FOM just as much as any other part of mathematics. A normative epistemological judgement, which is grounded solely on general philosophical principles and does not accurately take into account what falls under its scope (in our case this is FOM as they are practiced) can hardly be justified. 

The adequacy of the received Hilbert-style FOM is questionable on both these accounts. The most obvious problem concerns its \emph{feasibility}: it turns out that attempts to represent standard mathematical contents  in symbolic languages reserved by this FOM for this purpose require manipulations with chains of symbols, which cannot be represented in printed form and  surveyed by a human. A notorious example is Bourbaki's definition of cardinal number 1 that comprises about $10^{15}$ symbols \cite{Grimm:2010}. At first glance the unfeasibility problem appears to be merely ``practical'', and can be successfully dealt with by using symbolic shortcuts such as ``1'' . 
What is problematic here from an epistemological point of view is the relationships between mathematical theories represented ``informally'' or ``semi-formally'' with colloquial means  and their formal counterparts represented with formal languages that make part of Hilbert-style axiomatic FOMs.  The situation when theory $T$ is represented formally as $\lfloor T \rfloor$ , and the situation when $T$ is represented by different means but can   \emph{in principle} be translated into form  $\lfloor T \rfloor$, are very different from an epistemological point of view. In the latter case the given FOM have only a weak and indirect control over $T$, which works counter-positively. If $T$ is shown to be \underline{not} representable in the form  $\lfloor T \rfloor$ then one concludes that the given FOM does \underline{not} justify $T$. But the mere fact that $T$ \underline{is} representable in the form $\lfloor T \rfloor$ does not, by itself, allow one to judge whether $T$ is valid by the given FOM's standards or not. Such a judgement is possible only when the formal representation $\lfloor T \rfloor$ is built effectively.  Since  $\lfloor T \rfloor$ is usually unfeasible such a judgement is usually not available. This feature of the received Hilbert-style axiomatic FOM is at odds with the traditional expectation that mathematical theories can be justified by the very fact of being ``based on'' an appropriate FOM.  Justification of a given mathematical theory is a basic FOM-related task, which is, generally, not performed via metamathematical considerations concerning the epistemic (in-)completeness of the given theory.  Thus the feasibility problem has an epistemological dimension, which doesn't allow one to call this problem ``merely practical''. 


Another problem concerns basic epistemological assumptions about mathematics that underpin the received axiomatic conception of FOM. Hilbert assumes that (ready-made) mathematical knowledge is representable in the form of axiomatic theories, and that such theories are systems of sentences, where each sentence is deducible from a set of distinguished sentences called axioms. This conception assumes that all established mathematical knowledge has a propositional character. This strong epistemic assumption about the mathematical knowledge typically extends to all scientific knowledge and is not specific for Hilbert's axiomatic approach to FOM. As Bertrand Russell remarked in 1900 ``That all sound philosophy should begin with an analysis of propositions, is a truth too evident, perhaps, to demand a proof'' and followed this recipe in his seminal book \emph{The Principles of Mathematics} appeared in 1903 (for references and further discussion see \cite[p. v]{Gaskin:2009}). 

However the assumption according to which all mathematical knowledge deserving the name has a propositional character does not square well with the existing mathematical practice. This can be clearly seen with the historical example of Euclid's theory of elementary geometry presented in Books 1-4 of his \emph{Elements}. The first principles of Euclid's theory are of two types: geometrical \emph{postulates} and \emph{common notions} aka ``axioms'', which apply in the geometrical and in the arithmetical Books of the \emph{Elements}. Euclid's Postulates 1-3 are not propositions
\footnote{I assume here that a proposition, generally, admits for a truth-value. Euclid's Postulates 1-3 are not propositions because they don't admit for truth-values.} but schemes of elementary geometrical operations that intake some geometrical objects and output some other objects. For example, Postulate 1 is a scheme of operation that intakes two different points $A,B$ and outputs a straight segment  $AB$. Unlike the Postulates the Common Notions 1-3 do operate with propositions, namely with propositions of the form $A=B$ (equalities). However as Lewis Carrol rightly stressed back in 1895, these Common Notions should be read, on the pain of infinite regress, not as propositions but rather as schematic rules of inference, i.e., schemes of operations that allow one to obtain new equalities from some given equalities. Under this reading Common Notion 1 intakes two equalities of the form $A=C$ and $B=C$ and outputs equality of the form $A=B$. 

Accordingly, the main body of Euclid's geometrical theory comprises items called \emph{problems} along with theorems. Euclid's problems are recipes for making geometrical constructions provided with proofs that these constructions have their desired properties; these recipes are not expressed in a propositional form. For example, Problem 1 that opens Book 1 of the \emph{Elements} tells us how to construct an equilateral triangle given its side; it provides a construction using Postulates 1,3, and a proof (using Common Notion 1) that the output of this construction qualifies as an equilateral triangle. Most of Euclid's geometrical theorems also involve constructions, which in this case are colloquially called ``auxiliary''. So problems and theorems in Euclid's \emph{Elements} are intertwined and form a single deductive order, where earlier solved problems and earlier established theorems are used for solving further problems and establishing further theorems. Prima facie there is no reason for thinking that the operational \emph{knowledge-how} expressed in Euclid's \emph{Elements} in the form of problems is epistemically less significant than the propositional \emph{knowledge-that} expressed in Euclid's geometrical theorems\footnote{Logical reconstruction of Euclid-style geometrical reasoning remains problematic and cannot be more systematically explored here For a more detailed analysis of Euclid's geometrical reasoning see  \cite[Ch. 1.2]{Rodin:2014}. Back in 1974 Ian M\"uller claimed that no contemporary system of modern  logic could adequately account for Euclid's form of geometrical reasoning  \cite{Mueller:1974}, see also \cite{Mueller:1969}. In my opinion this situation did not significantly change during the following decades. }. 

Hilbert's work in the axiomatic foundations of Euclidean geometry \cite{Hilbert:1899} demonstrated that a rigorous propositional rendering of Euclid-style traditional geometry is possible. However it also showed that such a rendering is not innocent from a logical point of view, and requires a fundamental revision of Euclid's style of reasoning. In the Introductory part of their \emph{Foundations of Mathematics} published in 1934  Hilbert and Bernays stress the difference between their propositional axiomatic approach in theory-building and Euclid's ``genetic'' aka ``constructive'' approach, which involves non-propositional operations such as construction of new geometrical objects from given objects. The authors do not dismiss the genetic approach as an outdated version of their axiomatic approach but recognise its specific character, and conceive of a possibility to combine the two approaches \cite[v.1, Introduction]{Hilbert&Bernays:1934-1939}. 

The traditional Euclidean pattern of geometrical reasoning presented above remains very much alive in today's mathematics; various patterns of ``genetic'' reasoning co-exist in today's mathematics with the axiomatic reasoning  \cite[p. 35 ff]{Rodin:2014}. As we shall shortly (section 4.3 below) the Univalent Foundations of mathematics proposed in the mid-2000s provides a formal framework where the two patterns of mathematical reasoning work jointly.  Thus the assumption according to which all ready-made mathematical knowledge is essentially propositional, is at odds with the current mathematical practice where various forms of knowing-how continue to play a major role.

Let us now consider the normative dimension of this assumption.  A theoretical epistemological objection to this assumption can be formulated as follows. \emph{Models} are at least as important constituents of mathematical theories as the sentences that these models make true, and arguably are even more important (as will be explained shortly).  Today's mathematics is model-based with very little exception. Mathematical models (but not only mathematical theories construed as systems of sentences) play a key role in today's science. The thesis that mathematical knowledge reduces to propositional knowledge does not do justice to the epistemic role of mathematics in science and technology. 


Back in the 1950-ies Patrick Suppes put forward an influential epistemological  thesis according to which scientific theories should be identified with classes of models rather than with certain syntactic structures (like $\lfloor T \rfloor$) that are interpreted as systems of sentences. Since one and the same class of models can, generally, be determined by different theories identified syntactically, theories in Suppes' sense are invariant or ``robust'' with respect to a range of possible syntactic choices used for axiomatic presentation of these theories \cite{Suppes:2002}. Suppes' view of theories is known in the literature under the name of \emph{semantic} aka \emph{non-statement} view \cite{Halvorson:2016}. Suppes and his followers developed this view of theories attempting to apply Bourbaki-style axiomatic techniques in formal presentations of physical and other scientific theories. In spite of the fact that they referred to Bourbaki's project explicitly \cite[p. 33]{Suppes:2002}, \cite{Stegmuller:1979} the discussion of the concept of theory for some reason remained limited to philosophy of science and never extended to philosophical logic and philosophy of mathematics.  Similar views in philosophy of mathematics, however, were in different words expressed by proponents of \emph{mathematical structuralism} who were equally motivated by Bourbaki's project but preferred to use Bourbaki's term ``structure'' instead of the term ``model'' avoiding in this way a terminological interference with model theorists \cite{Hellman:2006a}.


The standard semantic approach in FOM doesn't do full justice to the model-character of mathematical reasoning because it offers a new way of thinking about existing formal techniques  (to wit, about the Hilbert-style axiomatic method empowered with Tarkski-style formal semantics) without developing new techniques appropriate for the purpose. True, building models (structures) and experimenting with them is a key element of today's mathematics. But the Hilbert-style axiomatic FOMs provide no means for the formal encoding such operations and giving them a stable theoretical form. This is because such FOMs comprise only formal rules, which are logical rules including rules of truth-preserving logical inferences.  FOMs of this type do not include rules similar to Euclid's Postulates that regulate manipulations with non-propositional objects.  Hilbert-style FOM may justify the \emph{existence} of a given mathematical structure with certain properties but is indifferent as to how this structure is obtained or ``found''. 
The idea that one may stipulate an object or structure with certain desired properties by fiat without asking how to build such an object from some familiar elements makes the axiomatic method in mathematics very useful and powerful. But the constructive
\footnote{Term ``constructive'' in logic an foundations of mathematics is highly overloaded. Following Hilbert and Bernays, I call here a mathematical reasoning constructive when it involves non-propositional operations such as building complex mathematical objects (aka constructions) from some given simple constituents, see \cite{Rodin:2018} for further discussion.} aka \emph{genetic} reasoning that does take such issues into account is at least equally important in mathematics. This constructive aspect of mathematics is particularly important in the applied mathematics where mathematical methods (applied outside of  pure mathematics) have a major epistemic role. As we shall see in section 4.3 below the Univalent FOM provides a technique for encoding mathematical knowledge-how and implementing this knowledge computationally.  


\section{Foundations at Work}
In this Section I provide a very brief overview of three different ``practical'' FOMs that emerged in the 20th and early 2st centuries: Nicolas Bourbaki's \emph{Elements of Mathematics}, the category-theoretic foundations first proposed by William Lawvere and the Univalent Foundations more recently proposed by Vladimir Voevodsky. This list is not supposed to be complete: there is a number of other interesting FOM-related developments in the 20th century mathematical practice that I do not consider here
\footnote{
This concerns, in particular, early attempts to adopt the Hilbert-style axiomatic presentation of elementary geometry in school education \cite{Halsted:1904} and current attempts to introduce computer-assisted proofs into mathematical practice, which are independent of Univalent Foundations \cite{Avigad:2018}.   
}.   
It is essential to keep in mind that we talk here not about some accomplished foundational theories per se but about continuing efforts to organise some large fragments of mathematics in a particular way with certain open-ended foundational theories: Set theory in the case of Bourbaki, Category theory in case of categorical foundations, and Homotopy Type theory (HoTT) in case of Univalent Foundations. The three projects started at different times: Bourbaki in the mid-1930s, Category-theoretic foundations in early 1960s and Univalent Foundations in the mid-2000s
\footnote{
More exact dates are, of course, conventional. Bourbaki group has been created in 1934 \cite{Baulieu:1994}. The idea of category-theoretic FOM has been first put forward by Lawvere in 1963 in his Ph.D. thesis \cite{Lawvere:1963}. The idea of Univalent Foundations (at this point without the name) has been first publicly presented by Voevodsky during his talk in the Institute of Advanced Studies (Princeton) on March 22, 2006 \cite{Voevodsky:2006}.  
}.  

Since the development of each later project involves a critique of earlier foundational projects and purports to fix their shortcomings, it is not unreasonable to regard these three projects as consecutive stages of the same overall process. However the real history of these projects does not fit into such a linear Hegelian scheme. Bourbaki-style set-theoretic FOM, category-theoretic FOM and Univalent Foundations (UF) all play roles in current mathematical practice \textemdash\ typically not in their ``pure'' forms but combined with some older foundational schemes and patterns. These roles are very different in the three cases. The set-theoretic FOM still maintain their ``official'' status (at least in eyes of the majority of logicians) and elements of Bourbaki-style presentation of mathematical reasoning are still widely used in practice (including mathematics education) in spite of the general scepticism about Bourbaki's foundational project shared by the majority of working mathematicians. Category theory is presently widely used as a convenient abstract theoretical language that extends and in some cases wholly replaces the Bourbaki-style set-theoretic language in many areas of mathematics. Situations when the Category theory has a role of \emph{foundation} are rarer and more specific; some such examples are referred to below. Finally, UF is a relatively new project (work in progress) aiming, in Voevodsky's words, at ``a computerised version of Bourbaki''.

\subsection{Bourbaki}
The continuing multi-volumed \emph{Elements of Mathematics} 
\footnote{The first volume of Bourbaki's \emph{Elements} was (first) published in 1939 \cite{Bourbaki:1939}. The latest original volume appeared in 2016 \cite{Bourbaki:2016}; in 2019 appeared a new edition of an earlier published volume. Further references and updates concerning the work of Bourbaki's seminar and recent publications of the group can be found on the group's webpage https://www.bourbaki.fr .}
produced  by a group of (mainly French) mathematicians of several generations that uses the pseudo-name of \emph{Nicolas Bourbaki} is a unique long-lasted attempt to build Hilbert-style axiomatic foundations of mathematics and implement them into a broad mathematical practice in research and education.  The title ``Elements of Mathematics'' refers Euclid's \emph{Elements}; following Euclid's example Bourbaki  aims at providing a self-contained  introduction into the contemporary mathematics, which systematically presents not only its principles but also its basic contents.  The original French title is \emph{El\'ements de math\'ematique}, which uses the unusual singular form `` math\'ematique'' that expresses Bourbaki's goal to unify mathematics.  A concise general description of this project, which makes explicit some grounding ideas behind it, has been published in 1950 as a separate programmatic article titled the ``Architecture of Mathematics'' \cite{Bourbaki:1950}
\footnote{The article has been written by Jean Dieudonn\'e and signed by the name of Nicolas Bourbaki, see \cite[ch.7]{Corry:2004}. }.

After a recognition of the unifying power of the Hilbert-style axiomatic method Bourbaki makes an interesting move by distinguishing between the \emph{logical} aspect of the axiomatic method from another aspect, which can be called \emph{structural}; in Bourbaki's view this is the latter rather than former aspect that makes the axiomatic method a powerful instrument of the unification: 

\begin{quote}
[E]very mathematical theory is a concatenation of propositions, each one derived from the preceding ones in conformity with the rules of a logical system  [$\dots$] It is therefore a meaningless truism to say that this ``deductive reasoning'' is a unifying principle for mathematics.  [$\dots$]  [I]t is the external form which the mathematician gives to his thought, the vehicle which makes it accessible to others, in short, the language suited to mathematicians; this is all, no further significance should be attached to it. What the axiomatic method sets as its essential aim, is exactly that which logical formalism by itself cannot supply, namely the profound intelligibility of mathematics.  \cite[p. 223]{Bourbaki:1950}
\end{quote}

In order to understand Bourbaki's misgivings about the ``logical formalism'' it is essential to take into account how Bourbaki's axiomatic method differs from Hilbert's. Using the term borrowed from the philosophy of science I have characterised in the last Section Bourbaki's method as \emph{semantic}. 
Bourbaki's ``axiomatic'' theories are \underline{not} Hilbert-style theories provided with some default interpretations but theories of set-theoretic models of such Hilbert-style theories. For example, each particular group (identified either up to set-theoretic equality or up to isomorphism) qualifies as a model of Bourbaki's axioms for Group theory (nowadays standard). But Group theory does not reduce to proving general theorems that follow from these axioms alone and thus apply to all groups indiscriminately. For example, the elementary Lagrange theorem does not belong to this kind because it involves two different groups, i.e., two models of the axioms. Since all such models are set-theoretic,  Lagrange theorem as well as any other group-theoretic theorem proved by Bourbaki is deducible from axioms of Set theory. This is how Bourbaki's semantic approach makes Set theory into a core of their FOM. 

Interestingly enough Bourbaki expresses some misgivings about Set theory along with the aforementioned misgivings about the ``logical formalism'':
 
\begin{quote}
[T]he axiomatic studies of the nineteenth and twentieth centuries have gradually replaced the initial pluralism of the mental [mathematical] representation  [$\dots$]  first to the concept of the natural number and then, in a second stage, to the notion of set. This latter concept, considered for a long time as ``primitive'' and ``undefinable'', has been the object of endless polemics [$\dots$]; the difficulties did not disappear until the notion of set itself disappeared [$\dots$] in the light of the recent work on logical formalism. From this new point of view, \emph{mathematical structures} [my emphasis - A.R.] become, properly speaking, the only ``objects'' of mathematics. \cite[p. 225-226]{Bourbaki:1950}
\end{quote}
 
The provided references do not allow one to identify the precise sense in which Bourbaki talks here about the ``disappearance of sets''.  But the above quotation makes it clear that already at an early stage of their project Bourbaki were motivated by the idea to represent mathematical structures without using Set theory. Yet, officially, all Bourbaki's structures are construed set-theoretically, namely, as set-theoretical models of certain systems of Hilbert-style axioms.  The intended notion of ``self-standing structure'' (without an underlying set) is nowhere rigorously defined in Bourbaki's \emph{Elements}. But as a matter of practical implementation of this intended notion Bourbaki systematically  disregard details of their own Set theory and its associated formal logical machinery.



A reason why Bourbaki feel so uneasy about their set-theoretic representation of mathematical structures concerns the idea (that historically dates back at least to Hilbert) that \emph{isomorphic} structures should count as identical \cite{MacLane:1996}. The problem is that no (Hilbert-style) axiomatic Set theory that comes with its internal identity relation provides a rigorous sense in which isomorphic structures can be said to be identical. We shall shortly see how this problem drives the category-theoretic approach in FOM and how it is settled in the Univalent FOM \cite{Ahrens&North:2019}.   In Bourbaki's  \emph{Elements} the notion of self-standing mathematical structure independent of any set-theoretic background remains  a wishful thinking or, on a more charitable interpretation, a \emph{regulative idea} in the Kantian sense of the word.

Outside a narrow group of enthusiasts who continue today to push further forward Bourbaki's project there are very few research mathematicians and mathematical educators who are ready to use Bourbaki's volumes in their teaching and build their research on them. Nevertheless it is clear that Bourbaki's \emph{Elements} indeed ``make explicit the essential general features, ingredients, and operations'' of certain style and pattern of mathematical thinking that has played a major role in 20th century mathematics. Even if a number of important achievements of 20th century mathematics have been wholly unrelated to Bourbaki's edifice, Bourbaki provided a certain vision of mathematics as a whole and contributed to its unification and connectedness by developing the semi-formal ``set-theoretic language'' that is still routinely used in various fields and areas of mathematical research and helps their cross-fertilisation. Even if there are serious reasons today to regard Bourbaki's set-theoretical FOM as outdated it still serves mathematicians and students of mathematics for many everyday purposes.

The story of the rise and fall of the Bourbaki-based approach in the elementary mathematics education, which took place worldwide in the course of 20th century, deserves special attention. It starts with an attempt at reforming school mathematical education on the national level according to the new standards known as the \emph{New Math} in the USA at the end of 1950s as a key element of the American response to the Soviet launch of the first Sputnik in 1957 \cite{Phillips:2015}. This educational reform aimed at the reinforcement of national technological and intellectual power and promoted the Bourbaki-style approach at the elementary school level.  It was led by a number of prominent American mathematicians including Marshall Stone and was predictably opposed by the majority of school teachers and pupils' parents. More remarkably, this reform was also severely criticised by some leading mathematicians and scientists, including such a prominent figure of the time as Richard Feynman \cite{Feynman:1965}. The opinion that the \emph{New Math} was a pedagogical error prevailed already in the early 1970s \cite{Kline:1973}, and this controversial reform was soon largely abandoned. 

Ironically or not,  in the Soviet Union a similar radical reform of mathematical education was started almost simultaneously with the American  \emph{New Math}.  In 1959 Boltyansky, Vilenkin and Yaglom published a programmatic paper \cite{Boltyansky&Vilenkin&Yaglom:1959} where they urged for renewal of standard school mathematical curriculum. In the late 1960s the reform of mathematical education was supported by Soviet authorities and soon developed into a project of national scale; a Bourbaki-style version of school mathematical curriculum was chosen for implementation in school textbooks, developed by a group of leading mathematicians and mathematical educators under the general supervision of Andrey N. Kolmogorov. The fate of this Soviet reform was similar to that of the American \emph{New Math}: in the late 1970s under the pressure of complaining school teachers, pupils' parents and harsh critique of some mathematicians, this reform was abandoned and the new mathematical textbooks were again replaced by more traditional ones.  

Similar reforms were undertaken in France and some other countries, and their fate was also similar. These national reforms of school mathematics did not proceed in isolation; they were parts of an international trend in mathematics education where the leading role belonged to French mathematicians including Andr\'e Lichnerowitcz and some other members of the Bourbaki group. 

The failure of Bourbaki-style approach in the elementary mathematics education cannot be fully analysed here. Let me only mention that proponents and critics of the Bourbaki-style approach in mathematics and mathematics education agreed that this approach was not helpful for pointing to, illuminating and supporting applications of mathematics in science and technology \cite{Stone:1961}, \cite{Feynman:1965}, \cite{Arnold:1998}. However they evaluated this fact differently: while the opponents argued that linking mathematics with science and technology was essential, the proponents of Bourbaki's approach argued that it was not.

\subsection{Category Theory as a Foundation}
Explaining the legacy of Set theory in the 20th century mathematics Yuri Manin  remarks that  in the second half of this century ``sets gave way to categories'' \cite[p.7]{Manin:2002}. The mathematical concept of \emph{category} emerged in 1940s, its official birth is marked by 1945 paper by S. Eilenberg and S. MacLane \cite{Eilenberg&MacLane:1945}. Category theory (CT) soon proved to be a useful mathematical language that helped to organise a number of emerging concepts and new results into a stable theoretical form. Among early examples that justified the usefulness of CT in eyes of many mathematicians were \cite{Eilenberg&Steenrod:1952}, \cite{Eilenberg&Cartan:1956}, \cite{Quillen:1967}. Some new mathematical concepts and mathematical theories developed during the same period would not even emerge without CT: in particular, this concerns the concept of topos invented by A. Grothendieck in 1950-ies.  A number of today's mathematical disciplines including Algebraic Geometry, Homological Algebra, Functional Analysis use CT as their basic language  \cite{Kromer:2007}, \cite{Marquis:2009}. 

The success of CT as an alternative language that allowed for effective theoretical accounting for mathematical structures in terms of their maps bypassing the  redundant set-theoretic details appeared to many mathematicians, among them some Bourbaki members (including H. Cartan, S. Eilenberg and A. Grothendieck) as a way to realise Bourbaki's dream about the ``disappearance of sets'' in their role of universal foundational background. Remarkably, Samuel Eilenberg and Norman Steenrod in the \emph{Preface} to their  \cite{Eilenberg&Steenrod:1952} as well as Daniel Quillen in \cite{Quillen:1967} boldly claim that using the language of CT they build ``axiomatic'' theories. Their approach was indeed axiomatic in its spirit in the sense that these authors stipulated the existence of categories with some special properties and on this basis developed abstract theories, which generalised upon many earlier known examples. However these theory did not meet standards of formal and logical rigour used by mathematicians working in FOM, and in fact did not pursue this goal.

The idea to develop CT into a genuine FOM is due to William Lawvere. In his 1963 Ph.D. thesis \cite{Lawvere:1963} he proposed Hilbert-style axioms for CT that are sometimes referred to in the current literature as Eilenberg-MacLane axioms and the corresponding theory is called EM; in the same thesis Lawvere built the Elementary Theory of Category of Sets (ETCS), where the word ``elementary'' refers to the first-order character of this theory, see also \cite{Lawvere:2005} for a more detailed exposition of this work. Unlike the aforementioned theories of Homological Algebra  and Homotopy \cite{Eilenberg&Steenrod:1952}, \cite{Quillen:1967} that also use the language of categories and are called by their authors ``axiomatic'', ETCS is an axiomatic first-order theory in the standard sense of the term familiar to logicians and philosophers. 
Later Lawvere developed a more general project of category-theoretic FOM, namely, an axiomatic theory of Category of Categories (Category of Categories as a Foundation or CCAF for short) \cite{Lawvere:1966a}. In these early works Lawvere challenged the notion according to which Set theory in some form should be a core part of modern FOM. In particular, he showed the possibility to develop an axiomatic Set theory on the basis of axiomatic CT rather than the other way round. 
However at least at this stage of his career he didn't challenge the notion according to which a modern FOM should involve the Hilbert-style axiomatic approach in some form. 


The most remarkable of Lawvere's achievement in the category-theoretic axiomatisation of mathematics is his work in Topos theory. The concept of \emph{topos} emerged in the circle of Alexandre Grothendieck in the early 1960s as a far-reaching generalisation of the concept of topological space. Grothendieck's concept of topos had no special relevance to logic; the discovery of logical structure of topos is wholly due to Lawvere and his foundational approach. In the beginning of his 1970 paper \cite{Lawvere:1970a} Lawvere provides his definition of topos usually called today the definition of \emph{elementary} topos; the title ``elementary'' reflects the fact that Lawvere's definition unlike Grothendieck's original construction almost straightforwardly translates into the standard first-order formal language \cite{McLarty:1992}.  A necessary and sufficient condition under which a given elementary topos is a Grothendieck topos was found in 1972 by Grothendieck's student Jean Giraud \cite{Giraud:1972}. 

Even if Lawvere's theory of elementary topos qualifies as a Hilbert-style axiomatic theory it has a feature, which makes it quite unlike ZFC and other similar theories. The standard approach applied by Hilbert and later described and justified by Tarski from a metodological viewpoint assumes that an axiomatic theory is based on a certain symbolic calculus with fixed logical semantics, which allows one to formulate some specific axioms and then interpret non-logical terms of this syntactic structure in certain intended models. The logical part of such theories, i.e., their ``underlying logic'' is supposed here to be rigidly fixed in advance \cite{Tarski:1941}. Lawvere's approach is different in this respect. His insight behind the axiomatisation of Topos theory is that every topos has an internal logical structure, which is analogous to that present in the category of sets.  As Lawvere puts it ``logic is a special case of geometry''. \cite[p. 329]{Lawvere:1970a}. This understanding of relationships between logic and geometry is very unlike Hilbert's and Tarski's. 


Studies in the category-theoretic FOM not only provided a new grounding for CT independent of Set theory but also helped to simplify many category-theoretic concepts and apply them in various areas of mathematics as well as outside of pure mathematics, in particular, in Computer Science.   Lawvere's axioms for elementary topos made Topos theory more accessible and  helped many people outside the community of experts in Algebraic Geometry to enter this field  and conduct fruitful research in it.  Although the category-theoretic FOM so far did not have a general impact on mathematics comparable with Bourbaki's impact it also had a significant effect on mathematics education (at the university level) and research.

\subsection{Univalent Foundations}

The Univalent Foundations (UF) is a novel foundation of mathematics first proposed by Vladimir Voevodsky in 2006 \cite{Voevodsky:2006}; UF uses Homotopy Type theory (HoTT) as a language and purports to express in this language various mathematical contents. HoTT in its original form is an interpretation of Martin-L\"of's Type theory (MLTT) in terms of Homotopy theory where types are interpreted as homotopy spaces and terms of these types are interpreted as points of those spaces. This interpretation was discovered by Vladimir Voevodsky and, independently, by Steven Awodey and Michael Warren in the mid-2000s \cite{UF:2013}.   

For Voevodsky a major motivation for developing UF was pragmatic and concerned the possibility of computer-assisted formal verification of long mathematical proofs. As we have already stressed, the received FOM does not support a formal checking and verification of mathematical proofs beyond trivial cases. However the practical need for such a reliable verification becomes more and more pressing. Even if the use of computer-assisted proofs remains today uncommon in the mainstream mathematical research, there is a growing corpus of  mathematical results proved with computer; in many such cases non-assisted ``purely human'' proofs are not known and there are indications that such proofs may not exist. A recent overview of existing proof-assistants and presentation of the state of the art in the interactive computer-assisted theorem proving \textemdash\ which however does not cover the  \emph{UniMath} and other UF-motivated projects \textemdash\ is found in \cite{Avigad:2018} and references therein. 

It is remarkable that no existing proof-assistants straightforwardly implement HIlbert-style formal representations of mathematical theories but use instead some form of Type theory in their theoretical background. Per Martin-L\"of designed and developed his constructive Type theory with dependent types (MLTT) back in 1980s already in view of its possible implementation on computer; its fragments were effectively implemented in Coq and some other proof-assistants. He provided this formal system with a \emph{proof-theoretic} logical semantics, which essentially differs from Tarskian model-theoretic logical semantics used in the received Hilbert-style FOM \cite{Piecha&Shroeder-Heister:2015}; syntactically, MLTT is a Gentzen-style (rule-based) but not Hilbert-style (axiom-based) formal system. The discovery of homotopical interpretation of MLTT in the mid-2000s, which gave rise to HoTT, convinced Voevodsky that this theory can serve as a basis for a new FOM, that supports an automated proof-checking. Voevodsky called this FOM \emph{univalent} because of the \emph{univalence axiom} (UA) explained in what follows. The origin of term ``univalence''  is somewhat arbitrary: it comes from the expression ``faithful functor'' translated into Russian on some occasion as (in the backward translation due to Voevodsky) ``univalent functor'' \cite[footnote 4]{Grayson:2018}. See \cite{Grayson:2018} for a concise introduction into UF, and \cite{UF:2013} for a more systematic exposition. 

The homotopical semantics of MLTT allows one to distinguish between types of different \emph{homotopy levels} and interpret types of different levels differently: types of level (-1) are interpreted as propositions, types of zero level are interpreted as sets, types of level 1 are interpreted as groupoids (namely, fundamental groupoids of homotopy spaces) and so on up to higher homotopy groupoids. Formal rules of the underlying Type theory  are interpreted, at the propositional level, as rules regulating logical inferences; at higher homotopy levels the same formal rules regulate geometrical (to wit homotopical) constructions, which make certain propositions true. Thus HoTT justifies  the aforementioned view of Lawvere on logic as ``a special case of geometry''. This feature of HoTT makes the UF-based mathematical reasoning in some respects similar to the traditional Euclid-style geometrical reasoning where geometrical constructions help to prove theorems. Since homotopy-theoretic concepts are spatial and intuitive in low dimensions, UF helps to visualise and understand many abstract mathematical concepts and proofs. In particular, it helps to have an intuitive grasp of formal proofs represented with a computer code and ready for an automated verification. This feature of UF addresses the issue of ``opaqueness'' (i.e., lack of surveyability) of computer-assisted proofs widely discussed in the literature \cite{Bassler:2006}.

Like in case of the traditional geometrical reasoning the epistemic role of UF-based constructions does not, generally, reduce to witnessing certain mathematical sentences but have other important epistemic roles. Notice that by constructions I mean here rule-based mathematical procedures but not only their outcomes. Such constructions represent mathematical know-how in an explicit form. For a historical example think of Euclid's construction of five regular polyhedra aka Platonic Solids in Problems 13.13-17 of \emph{Elements}, Book 13. This is the crowing achievement of Euclid's edifice; when one traces these constructions back to their first principles, i.e., to Postulates and Common Notions of Book 1 then one can see that the necessary theoretical background of these five Problems comprises a significant portion of the mathematical knowledge presented in the \emph{Elements}. This large bulk of mathematical knowledge is packed in the five aforementioned Problems in the operational but not propositional form.  For another important historical example think of algebra in its early form presented by al-Khwarizmi in the book, which gave the discipline its current name \cite{Rashed:1994}. Al-Khwarizmi's algebra is a systematised collection of recipes  (provided with an Euclide-style geometrical background) for solving problems of a certain sort. The current mathematical practice does not give one any indication that the role of the operational knowledge-how in mathematics is less significant today than it used to be in the past. The fact that UF unlike the standard Hilbert-style axiomatic FOMs allow for a formal and meticulous representation of operational knowledge in the form of UF-based homotopical constructions is one more reason in favour of the claim that UF are more adequate.


The significance of knowing-how is particularly evident in the applied mathematics where mathematics is seen as an epistemic tool rather than an independent discipline. Mathematical models of natural phenomena and technological artefacts always verify certain mathematical sentences;  to put it colloquially, such models have certain mathematical properties. However the main epistemic role of such models is not to instantiate and witness these mathematical properties but to account for empirical data. In such contexts knowing \emph{how} to build an appropriate mathematical model is a key issue.  

Let me briefly illustrate this point with two examples: the historical example Ptolemean Astronomy and the up-to-date example of Topological Data Analysis (TDA). Ptolemean Astronomy involves the traditional Euclid-style geometry as its mathematical foundation (but goes beyond the geometrical content of Euclid's \emph{Elements} since it applies Spherical geometry). Visible motions of  celestial bodies are reconstructed as superpositions of circular motions. In case of planets this requires using a number of supplementary circular motions, which are not observable. The trajectories of such motions Ptolemy calls \emph{epicycles} \cite{Neugebauer:2012}. The resulting geometrical model, by Ptolemy's famous word, ``saves the phenomena'', i.e., is empirically adequate (with a certain degree of approximation). Importantly, Ptolemy's geometrical models of astronomical phenomena are built step by step and thus adjusted to observational data; such genetic procedures are based on constructive geometrical principles closely related to Euclid's Postulates (think of epicycles).  Such a systematic way of model-building is made possible by the constructive (genetic) character of Euclid-style geometry explained in section 3 above. 

In TDA, which is a rapidly developing area of today's applied mathematical research \cite{Patania&Vaccarino&Petri:2017}, large datasets of various nature are accounted for in terms of topological spaces and their invariants. This allows one to extract from the given data some essential information, which turns out to be relevant and useful in many cases. All relevant characteristics of those topological spaces \textemdash\ and, in a sense, also these very  spaces \textemdash\ are computed on the basis of the given data. Even if such a procedure can, in principle, be fully specified in terms of truth-conditions within an appropriate external theoretical framework, UF offers a possibility to interpret such computations internally (and yet formally and rigorously) as genetic topological  constructions. Such a direct internal interpretation better squares with the colloquial language and intuition of working mathematicians and computer scientists working in the field. This specific feature of UF has a potential to narrow the existing gap between today's abstract pure mathematics and applied mathematics \cite{Rodin:2021}.

Let us now turn to the Univalence Axiom (UA), which gives the Univalent Foundations its name. Let $A, B$ be types, in symbols $A, B \: TYPE$. Consider the identity type $A=_{TYPE}B$ and the type of \emph{equivalences} $A \simeq_{TYPE} B$ where by equivalence one understands  a function $f: A \rightarrow B$, which is in an appropriate sense invertible; under the homotopical interpretation such equivalences are homotopy equivalences. The rules of MLTT/HoTT allow one to construct a canonical map of the form

$$e \: (A = B) \rightarrow (A \simeq B)$$      

\noindent which witnesses the fact that identity is a special case of equivalence. The \emph{Univalence Axiom} states that this map $e$ has an inverse and thus is itself an equivalence. In other words, UA says that the type
 
$$(A = B) \simeq (A \simeq B)$$

\noindent is inhabited. More recent versions of UF based on Cubical Type theory (CTT) do not use UA as an axiom but prove it as a theorem.  The aforementioned inverse map, which is not constructible by using the rules of MLTT, is constructible by using the rules of CTT. Thus in CTT-based UF there is no need to postulate univalence with a special axiom. This fact supports the claim that CTT is a better formal carrier for UF that MLTT \cite{Cohen&Coquand&Huber&Mortberg:2016}. 

Univalence is important for a variety of reasons. A  reason which has a particularly strong philosophical appeal is that the univalence property apparently justifies the  structuralist idea of mathematical reasoning ``up to isomorphism'' (or up to another type of equivalence), which reflects the existing practice in the Bourbaki-style set-based mathematics and in the more recent  category-theoretic mathematics but lacks a firm logical grounding in these foundational frameworks.  A detailed analysis of this claim, which includes further qualifications and certain reservations, is found in \cite{Ahrens&North:2019}.

It should be born in mind that the homotopical semantics of MLTT is an extension (not fully conservative) of its intended logical semantics but is not a model of this theory in Homotopy theory. \emph{Models} of HoTT are further interpretations of its syntactic structures, which respect the homotopical semantics in an appropriate sense. When formal rules of HoTT are interpreted in external structures they can be understood as rules for building such structures from their elements. This feature of UF has no analogue in the standard Bourbaki-style semantic approach to theory-building where all structures are supposed to be ``static''. From a historical perspective this aspect of UF can be described as a synthesis of axiomatic  and \emph{genetic} approach, which has been considered by Hilbert and Bernays back in 1934 as a theoretical possibility \cite{Rodin:2018}.      



Since UF is still a work in progress it is too early to evaluate its success and practical impact on mathematics against the impacts of earlier FOMs. However already today one can see that UF is designed as a practical tool for learning and doing mathematics, which unlike the received Hilbert-style axiomatic FOMs does not suffer from the ``meta-mathematical bias'' (see section 2 above). Thus UF presently appears as a strong candidate for bridging the existing wide gap between FOM and the rest of mathematics, between mathematical foundational theories and mathematical practices.

\section{Renewing Foundations and Progress in Mathematics}

In the remaining part of this paper I would like to describe my general views about mathematics and its foundations in terms proposed by Michele Friend \cite{Friend:2014}.  First of all, my position is \emph{foundationalist} insofar as I believe that FOM is an important part of mathematics both practically and theoretically; I also believe that FOM is an important (albeit not unique) point of contact between mathematics and its philosophy. Second of all, I am a \emph{pluralist} about FOMs. This should be clear from the earlier sections where I described several different FOMs without saying which is the right one, and which are wrong. 

Since my working general concept of FOM borrowed from Lawvere and Rosenburgh is broader than the one used by Friend, my position is foundationalist in a weaker sense of this term than in Friend's sense.  Friend defines FOM as ``an axiomatically presented mathematical theory to which all or most of successful existing mathematics can be reduced'' \cite[ p. 8]{Friend:2014}. I assume that Friend uses here the term ``axiomatic'' in its usual Hilbertian sense, and I qualify an axiomatic presentation of mathematics in \emph{this} sense as a special feature of the received Hilbert-style FOMs rather than as an essential feature of FOM in general. Univalent Foundation is an example of FOM that is not axiomatic in this usual sense because the UF applies a Gentzen-style (i.e., rule-based) but not a Hilbert-style formal architecture
\footnote{It makes sense to broaden the received notion of being axiomatic, so it may include also Gentzen-style formal systems. This terminological convention, however, creates a high risk of confusion (since systems like MLTT that comprise no axiom in the received sense of this term would also qualify as axiomatic). For this reason I stick in this paper to the received meaning of ``axiomatic''. For discussion of a broader conception of being axiomatic see my \cite{Rodin:2018}}. 

Further, I don't assume together with Friend that mathematics can be ``reduced'' to its foundations in any epistemically relevant sense of this term. When such a reduction is understood in a strong epistemic sense, so the reducibility condition reads

\emph{In order to learn all existing mathematical knowledge it is sufficient to learn FOM.} 

this condition is obviously too strong and even absurd. When the reduction is understood as formal provability or in some other sense that has no definite epistemic content, the same condition is too week or irrelevant. The fact that all texts written in English can be encoded into 0-1 sequences does not make 0 and 1 into building blocks of all English literature. Similarly, the fact that a lot of mathematics can be encoded into fragments of ZFC does not, by itself, make ZFC into a foundation of mathematics
\footnote{I'm not trying here to dismiss the set-theoretic FOM in one stroke and realise that set-theoretic FOMs involve more than formal provability. I would like to stress, however, that formal provability of theorem $T$ in ZFC in the usual sense of existence of truth-preserving inference from the axioms of ZFC to $T$ does not, by itself, shed a light on how $T$ is (or is not) \emph{known}, see \cite{Prawitz:1986} for a discussion. }
Instead of looking further for an appropriate notion of reduction I shall describe the relationships between mathematics and its foundations by looking at it from the opposite end. The following epistemological condition on FOM is, in my view, appropriate: 

\emph{In order to acquire any fragment of the existing mathematical knowledge it is \emph{necessary} to learn (a relevant fragment of) FOM .}

\noindent The latter condition implies that FOM or its relevant fragment are necessary in justification of every mathematical result but it does not imply that FOM are sufficient for that purpose.   

In spite of these differences between my working concept of FOM and Friend's concept, they share many essential features in common, in particular, a normative character, a capacity to unify mathematics, and a justificatory function. So my notion of being a foundationalist is not wholly at odds with Friend's, after all. 

Let me now explain how I think about multiple FOMs. My first point is historical. I believe that the development of FOM is an important part of the overall historical development of mathematics. I assume that this overall development is progressive, and understand this progress as a growth of accumulated mathematical knowledge \cite{Bird:2007}.This notion of progress is compatible with a picture of mathematics that develops on some fixed foundations (say, Euclid's), so the growth of mathematical knowledge amounts to proving more and more geometrical theorems and solving more and more geometrical problems by ruler and compass. This process is potentially infinite, at least if we abstract away from issues of feasibility. However mathematics as we know it is more interesting, it simply does not develop in this simplistic way, at least not in the long term. Our real mathematics involves a permanent \emph{revision} of foundations, which contributes to the overall mathematical progress as follows.  

Recall that the trisection of angle and some other geometrical problems are unsolvable on Euclid's foundations. These old geometrical problems were successfully solved (in the negative way) in the 19th century using the then new algebraic foundations of geometry laid out during the course of 17-18th centuries by Descartes and others. These new foundations \cite{Arnauld:1683} were not a conservative extension of Euclid's foundations with certain additional principles but involved a very deep conceptual revision of these old foundations. However this revision did not disqualify Euclid's foundations as being wholly erroneous either. Euclid is still with us, and we see nothing fundamentally wrong  about his geometrical constructions by ruler and compass (except for some fixable logical gaps). This situation can be compared with the foundational change in Physics that was brought about with the development of Quantum Mechanics: this foundational theory did not replace the earlier accepted foundational theory (Classical Mechanics) but explained the success of this older theory while showing the limits of its applicability. The principles according to which new foundational theories in physics should relate to older foundational theories in this way is known as Bohr's  \emph{correspondence principle}. The same principle applies to foundational changes in mathematics.   

In order to use the new foundations for solving geometrical problems formulated in a different foundational framework these problems had to be reformulated in a new way.
What we learn today under the name of Pythagorean Theorem is not quite the same statement that we typically identify as this same theorem in Euclid's \emph{Elements} (Book 2, Theorem 47); the corresponding proofs of these statements are moreover very different. How and in which exact sense mathematical contents preserve their identities through such foundational changes, is a tricky and little explored question, which I should leave aside. However the very fact that such a preservation of mathematical knowledge survives foundational changes is remarkable. 

The above example of foundational change in geometry is typical rather than exceptional. Sometimes new foundations help to identify old errors 
\footnote{This is the case of Hilbert's foundations of geometry, which provided a new general proof that all earlier proposed alleged proofs of Euclid's Fifth Postulate were erroneous. In fact, however, those alleged proofs were found to be erroneous or at least suspicious and unreliable even before this foundational change.},
but I don't know about examples of foundational revisions that have invalidated large bodies of mathematical knowledge established earlier on different foundations. Thus I observe that a revision of FOM typically has these three effects: (i) it preserves earlier acquired mathematical knowledge but expresses it in a new form, (ii) it helps to solve open problems and (iii) it triggers new developments that were not possible with only the resources of the old foundations. So the continuing revision of FOM contributes to the overall progress of mathematics.   

The above historical point does not commit me to a historical relativism about epistemic norms provided by different FOMs used in different historical times. It is quite justified to judge older FOMs by epistemic standards of more recent FOMs. One should bear in mind, however, that different FOMs often have very different designs, so some such judgements can be easily ill-formed. To repeat, new FOM, generally, does not disqualify the older FOM that they replace but put them in a new perspective. I can see no reason for thinking such a renewal of FOM should ultimately converge to a certain ``true'' FOM. The continuing revision of FOMs has a character of expansion (in its depths and scope) rather than convergence.  I believe that this specific form of progress is worth encouraging rather than being stopped by choosing the current FOM as ultimate.

\section{Tolerance and Pluralism}

My final point concerns what Friend calls the principe of \emph{tolerance} to conflicting philosophical views on mathematics associated with different FOMs \cite{Friend:2014}. What has been said above about the continuing revision of FOMs during the historical development of mathematics should not be understood in the sense that at any given moment of historical time mathematics relies or should rely on a unique FOM. In practice,  mathematicians of the same generation often use different foundations, including more traditional and more ``progressive'' ones.  The failure of the New Maths and akin reforms in mathematics education in the 1960-1970s suggests that such a diversity is justified at least in education. I have described above how new foundations replace and ``absorb'' older foundations without disqualifying these older foundations completely. A similar scheme applies to contemporary FOMs. There exists a large body of mathematical knowledge that can be reconstructed on different alternative FOMs, say, on the set-theoretic and the category-theoretic FOMs. Let us consider a more problematic case when certain mathematical proof (construction) $P$ is justified by one FOM ($F$) and not justified by another FOM ($G$). It appears that such cases make a tolerant attitude to conflicting FOMs impossible on pain of contradiction
\footnote{Tolerating contradictions by using a paraconsistent logic is an option but I'm going to propose here a different solution.}
: in order to judge whether or not $P$ is valid one need to take a side and reject either $F$ or $G$.  This contradiction, however, can be avoided. There is a powerful mechanism that allows one to avoid the choice between $F$ and $G$ without becoming an epistemological relativist or a sceptic.            

A more precise description of the above conflict between $F$ and $G$ involves  two different versions of $P$, namely, $P_{F}$ that construes and justifies $P$ according to standards of $F$, and $P_{G}$ which construes $P$ according to standards of $G$ and invalidates it. Now suppose that  $G$ admits an interpretation $g$ of $F$ that extends to $P_{F}$ and provides an explanation in terms of $G$ of how and why proof $P_{F}$ is valid by the standards of $F$. Reciprocally, $F$ admits interpretation $f$ of $G$ that extends to $P_{G}$ and explains in terms of $F$ why $P_{G}$ is invalid by the standards of $G$. When such interpretations $f, g$ exist $F, G$ are called \emph{bi-interpretable}. These mutual interpretations can be pictured with the following diagram:

$$\xymatrix{F/P_{F}\ar@<2pt>[r]^{g} &G/P_{G}\ar@<2pt>[l]^{f}}$$

\noindent So a user of $G$ gets another reconstruction of $P$ in terms of $G$, namely, $P'_{G} = g(P_{F})$ which unlike $P_{G}$ is valid by the standards of $G$. In other words, the interpretation $f(F)$ makes it possible to fix $P_{G}$, which provides a reconstruction of $P$, which is valid in terms of $G$. Reciprocally, a user of $F$ gets another reconstruction of $P$, namely, $P'_{F} = f(P_{G})$ which unlike $P_{F}$ is invalid by the standards of $F$. The users of $F$ and $G$ now share a mutual understanding of principles needed for justifying or rejecting $P$. The user of $F$ and the user $G$ will think about $P$ in different ways and disagree about the epistemic status of $P$.  But at this point their disagreements belong to the philosophy of mathematics rather than mathematics itself. 

Let $P$ be proof of proposition (theorem) $T$, and $T_{F}, T_{G}$ be reconstructions of $T$ in languages of $F$ and $G$ correspondingly. Since $P_{F}$ is valid by the standards of $F$ a user of $F$ is justified to believe that $T_{F}$ is true. For the sake of argument let us further assume that $T_{G}$ is provably false, so a user of $G$ is justified to believe it. As far as  $T_{F}$ and $T_{G}$ are seen as two different representations of the same proposition $T$ it can be argued that users of $F$ and users $G$ have beliefs, which contradict one another, so a person who would wish to use both these FOMs simultaneously would have contradictory beliefs. However this reasoning is erroneous. Propositions $T_{F}$ and $T_{G}$ are not the same; a self-standing mathematical proposition $T$, which is given without a theoretical and foundational context, is ambiguous and has no definite true-value. (Think of such mathematical terms as ``space'', ``point'', ``circle'', ``set'', which outside a theoretical context have no definite meaning.) Rendering $T$ as $T_{F}$ and, alternatively, as $T_{G}$ disambiguates $T$ by splitting it into two different propositions. The above scheme allows one to construe this disambiguation internally for $F$ and for $G$ : a user of $F$ construes it as the distinction between $T_{F}$ and $f(T_{G})$ while a user of $G$ construes it as the distinction between $T_{G}$ and  $g(T_{F})$. This helps us to avoid pointless arguments about which of $T_{F}, T_{G}$ is the ``right'' rendering of $T$. Modulo the equivalences $T_{F} \simeq g(T_{F})$ and  $T_{G} \simeq f(T_{G})$, which hold when interpretations $f, g$ are faithful, users of $F$ and $G$ agree that  $T_{F}$ is true and $T_{G}$ is false.        

The above scheme does not apply universally but it applies to bi-interpretable FOMs, which include many cases of interest. For a trivial example consider the case when $F, G$ are two different Hilbert-style axiomatic theories which share some axioms in common. One can use ZF as a foundation and use the Axiom of Choice (AC) in certain proofs as a separate hypothesis. Alternatively, one can take ZFC as a foundation, and this time view the AC as an axiom. From an epistemological and foundational point of view the two approaches differ. There is, however, a large body of mathematical knowledge that is wholly neutral with respect to choice between accepting AC as an axiom and using it as a hypothesis. In this case the interpretations $ZF \rightarrow ZFC$ and $ZF \leftarrow ZFC$ are trivial and amount to changing the epistemic status of AC from hypothesis to axiom and back. In cases when different FOMs use different formal architectures (as in the case of ZF and UF) and differ conceptually such mutual interpretations can be very involved and non-trivial. The long-lasted continuing history of rebuilding Euclidean geometry on different foundations provides a lot of suggestive examples of this sort.

Some conditions stronger than the mere bi-interpretability can also be considered. Let $f,g$ be composable in a associative way.  If $g\circ f = 1_{F}$ and $g\circ f = 1_{G}$, i.e., the image of image of $F$ in $G$ is $F$ itself and, reciprocally, the image of image of $G$ in $F$ is $G$ itself, the two FOMs are equivalent \textemdash\  in a way, which is specified by the character of these interpretations. The elementary Category theory, this time in its role of a meta-theoretical tool, suggests many other possible configurations between different FOMs, which I leave here for another study. 

The possibility of the co-existence of different FOMs in mathematical practice does not undermine the significance of philosophical arguments in favour of certain FOMs and against certain other FOMs. Such philosophical discussions have multiple impacts on mathematics including shaping its new trends. The tolerant attitude to multiple FOMs and their philosophical underpinnings, which I defend here,  does not imply a neutrality and lack of interest to FOM-related philosophical matters. I understand tolerance as a methodological principle according to which a foundational debate in mathematics needs an environment, which is both concurrent and cooperative. As far as one accepts that conflicting philosophical ideas about mathematics deserve a common stage where they are systematically presented, criticised and further developed, one should also accept the notion that a part of this stage belongs to mathematics and accommodates multiple FOMs. As I have explained above, the multiplicity of FOMs does not break mathematics down into disconnected pieces. Different FOMs may serve different special purposes and be useful in different ways. We showed that the received set-theoretic axiomatic  FOM was designed for meta-theoretical purposes and that this FOM cannot be effectively used for the purpose of formal proof-verification. Using UF for this latter purpose does not exclude a possibility to use the received set-theoretic FOM for meta-mathematical purposes. I don't believe that these different epistemic tasks \textemdash\ to verify proofs and to provide a meta-theoretical view on a given theory \textemdash\  can be adequately distinguished as being practical and being theoretical. Both tasks have theoretical aspects that are analysed in normative logical and epistemological terms and practical aspects that are described in terms of past experiences and recommendations for novices.

As many authors have  earlier remarked, the traditional architectural metaphor behind the concept of theoretical foundation can be not only useful and suggestive but also very misleading. Unlike traditional buildings mathematical theories do not rise and fall with their foundations but have a capacity to survive through foundational changes. Martin Krieger compared mathematics of the 19th century with an industrial city \cite[Ch. 6]{Krieger:2015}. In line of Krieger's metaphor FOM should be thought of as a city infrastructure, which uses elements inherited from earlier historical times. We may extend Krieger's metaphor and think about today's FOMs as even larger, more flexible and more robust infrastructures, which are more suitable for our diverse post-industrial economies. Pluralism in Friend's sense is a \emph{sine qua non} of FOM so conceived. My short answer to the question asked in the title of this paper is this: one mathematic(s) and multiple foundations.   
      
\bibliographystyle{plain} 
\bibliography{fom20}

\end{document}